\documentclass[12pt]{article}
\usepackage{amsfonts}
\usepackage{}
\usepackage{amsthm,amsmath}
\usepackage{amscd}
 \usepackage{amssymb}
 \usepackage{hyperref}
 \usepackage[all]{xypic}
 \usepackage{mathtools}
  \usepackage[numbers,sort&compress]{natbib}
 \usepackage[utf8]{inputenc}
  \newtheorem{lemma}{Lemma}[section]
 \newtheorem{corollary}[lemma]{Corollary}
 
 \newtheorem{theorem}[lemma]{Theorem}

\usepackage[utf8]{inputenc} 

\def\includegraphics{}

\usepackage{tikz,xcolor,hyperref}

\definecolor{lime}{HTML}{A6CE39}
\DeclareRobustCommand{\orcidicon}{
\begin{tikzpicture}
\draw[lime, fill=lime] (0,0)
circle[radius=0.16]
node[white]{{\fontfamily{qag}\selectfont \tiny \.{I}D}};
\end{tikzpicture}
\hspace{-2mm}
}
\foreach \x in {A, ..., Z}{%
\expandafter\xdef\csname orcid\x\endcsname{\noexpand\href{https://orcid.org/\csname orcidauthor\x\endcsname}{\noexpand\orcidicon}}
}

\usepackage{amscd}
\usepackage{amsmath,amssymb,mathrsfs}
\usepackage{hyperref}
\usepackage[all]{xypic}
\usepackage{times}
\usepackage{inputenc}
\usepackage{titlesec}

\allowdisplaybreaks[1]

 \usepackage{amsthm}

\advance\headheight1.15pt
\newenvironment{proof of Theorem 1.1 step1}{{\noindent\it {\bf Proof of the implication $(1)\Rightarrow(3)$}}\quad}{\hfill $\square$\par}
\newenvironment{proof of Theorem 1.1 step2}{{\noindent\it {\bf Proof of the implication $(3)\Rightarrow(2)$}}\quad}{\hfill $\square$\par}
\newenvironment{proof of Theorem 1.2}{{\noindent\it {\bf Proof of Theorem 1.2}}\quad}{\hfill $\square$\par}
\newcommand{\sn}{\sum_{n=0}^{\infty}}
          
\newcommand{\dd}{\mathbb{D}}   \newcommand{\ba}{\mathcal {B}^{\alpha}}   \newcommand{\bb}{\mathcal {B}^{\beta}}
\newcommand{\Cu}{\mathcal{C}_\mu}      
\newcommand{\ii}{\int_{0}^{1}}    \newcommand{\hd}{H(\mathbb{D})}
\newcommand{\comment}[1]{}

\begin{document}
\baselineskip=8pt

\title{\fontsize{15}{0}\selectfont  Ces\`aro-like operator acting between Bloch type  spaces}

\author{\fontsize{11}{0}\selectfont
   Pengcheng Tang$^{*,a}$ and Xuejun Zhang$^{b}$
   \\
   \fontsize{10}{0}\it{
 $^{a}$School of Mathematics and Statistics, Hunan University of Science and Technology, Xiangtan, Hunan 411201, China}
 \\
   \fontsize{10}{0}\it{$^{b}$ School of Mathematics and Statistics, Hunan Normal University, Changsha, Hunan  410081, China}
}

\date{}
\maketitle
\thispagestyle{empty}
\begin{center}
\textbf{\underline{ABSTRACT}}
\end{center}
 \ \ \ \ \ \ Let $\mu$ be a  finite positive Borel measure on the interval $[0,1)$ and $f(z)=\sum_{n=0}^{\infty}a_{n}z^{n} \in H(\mathbb{D})$. The Ces\`aro-like operator is defined by
$$
\Cu (f)(z)=\sum^\infty_{n=0}\left(\mu_n\sum^n_{k=0}a_k\right)z^n, \ z\in \dd,
$$
where, for $n\geq 0$, $\mu_n$ denotes the $n$-th  moment of the measure  $\mu$, that is,
$\mu_n=\int_{[0, 1)} t^{n}d\mu(t)$.

 In this paper we investigate the action of
the operators $\Cu$  from one Bloch type  spaces $\ba$ into another one $\mathcal {B}^{\beta}$.

\begin{flushleft}
\end{flushleft}
\begin{flushleft}
\end{flushleft}

\let\thefootnote\relax\footnote{$^*$Corresponding Author}
\let\thefootnote\relax\footnote{ Pengcheng Tang: www.tang-tpc.com@foxmail.com}
\let\thefootnote\relax\footnote{ Xuejun Zhang: xuejunttt@263.net}

\vspace{1cm}

\section{Introduction} \label{Sec:Intro}

 \ \ \  \ \ \ Let $\mathbb{D}=\{z\in \mathbb{C}:\vert z\vert <1\}$ denote the open unit disk of the complex plane $\mathbb{C}$ and $H(\mathbb{D})$ denote the space of all analytic functions in $\mathbb{D}$.

 For $0<\alpha<\infty$, the Bloch-type space, denoted by $\ba$, is defined as
 $$\ba=\{f\in \hd:||f||_{\ba}=|f(0)|+\sup_{z\in\dd}(1-|z|^{2})^{\alpha}|f'(z)|<\infty\}.$$

When $\alpha=1$, $\ba$ is just the classic Bloch space $\mathcal {B}$.

For $0<\alpha<1$ , the analytic Lipschitz space $\Lambda_{\alpha}$ consists of  the functions  $f\in \hd$ for which
$$||f||_{\Lambda_{\alpha}}=\sup\left\{\frac{|f(z)-f(w)|}{|z-w|^{\alpha}}: z,w \in \dd, z \neq w\right\}<\infty.$$
It is known that  (see \cite{b1}) $\ba \cong\Lambda_{1-\alpha}$ for $0<\alpha<1$ and  $\Lambda_{\alpha}$ is  contained in the disc algebra.

 For $f(z)=\sum_{n=0}^\infty a_nz^n\in \hd$, the Ces\`{a}ro operator $\mathcal {C}$ is defined  by
 $$
\mathcal {C}(f)(z)=\sum_{n=0}^\infty\left(\frac{1}{n+1}\sum_{k=0}^n a_k\right)z^n, \  z\in\dd.
$$

 The boundedness of the Ces\`{a}ro operator has been studied by several authors on certain spaces
of  analytic functions. See, e.g.,\cite{ces5,ces10,ces3,ces12,ces13,ces15,ces4,ces14} and the references therein. The Ces\`{a}ro operator $\mathcal {C}$ has also been generalized to various forms and  its generalization has been widely studied on the space of holomorphic functions. For instance,   Stevi\'{c} \cite{st0}   studied  the generalized Ces\`{a}ro operator on the polydisc.  Hu \cite{ces11} studied the extended Ces\`{a}ro operators on  the Bloch space in the unit ball of $\mathbb{C}^{n}$.  Stevi\'{c}  also considered the generalized   Ces\`{a}ro operator on weighted-type spaces in \cite{st6} and studied the generalized Ces\`{a}ro operators acting between Bloch type spaces in \cite{st7}.    For more information on some generalizations of the Ces\`{a}ro operator on  spaces of holomorphic functions,  the reader is referred to \cite{ces9,ces8,ST1,ST2,ST3,st4,st5}.



 Recently,  Galanopoulos, Girela and Merch\'an  \cite{ces1}  introduced  a Ces\`aro-like operator $\Cu$ on $\hd$, which is a natural generalization of the classical Ces\`{a}ro operator $\mathcal {C}$. They  systemically studied the  operators $\Cu$  acting on distinct spaces of analytic functions, such as Hardy space, Bergman space, Bloch space.

  Let $\mu$ be a positive finite Borel measure on $[0, 1)$ and  $f(z)=\sum_{n=0}^\infty a_nz^n\in \hd$. The Ces\`aro-like operator $\Cu$ is defined  as follows:
  $$
\Cu (f)(z)=\sum^\infty_{n=0}\left(\mu_n\sum^n_{k=0}a_k\right)z^n=\int_{0}^{1}\frac{f(tz)}{1-tz}d\mu(t), \  z\in\dd.
$$
where $\mu_{n}$ denote the moment of order $n$ of $\mu$, that is,
  $\mu_{n}=\ii t^{n}d\mu(t)$.
If $\mu$ is the Lebesgue measure on $[0, 1)$, the operator $\Cu$ reduces to the classical Ces\`{a}ro
operator $\mathcal {C}$.


The Ces\`aro-like operator  $\Cu$ can be regarded as  an operator induced by the matrix
\[
	         \Cu=\begin{pmatrix}
	         \mu_{0}& 0 & 0& 0& \cdots\\
	       \mu_{1}& \mu_{1} & 0&  0& \cdots\\
           \mu_{2}& \mu_{2} & \mu_{2}&  0& \cdots\\
           \vdots& \vdots& \vdots& \vdots& \ddots
	         \end{pmatrix}.
	         \]
\[
	    \begin{pmatrix}
	         \mu_{0}& 0 & 0& 0& \cdots\\
	       \mu_{1}& \mu_{1} & 0&  0& \cdots\\
           \mu_{2}& \mu_{2} & \mu_{2}&  0& \cdots\\
           \vdots& \vdots& \vdots& \vdots& \ddots
	         \end{pmatrix}
     \begin{pmatrix}
     a_{0}\\a_{1}\\a_{2}\\ \vdots
           \end{pmatrix}
           =  \begin{pmatrix}
    \mu_{0}a_{0}\\ \mu_{1}\sum_{k=0}^{1}a_{k}\\ \mu_{2}\sum_{k=0}^{2}a_{k}  \\ \vdots
           \end{pmatrix}.
	         \]

The  Ces\`aro-like operator $\Cu$ defined above has  attracted the interest  of many mathematicians.
Jin and Tang \cite{ces6} studied  the boundedness and compactness of  $\Cu$  from one Dirichlet-type space $\mathcal {D}_{\alpha}$  into another one $\mathcal {D}_{\beta}$.
Bao, Sun and Wulan  \cite{baoo} studied the range of $\Cu$ acting on $H^{\infty}$.  Blasco \cite{blas,blas2} investigated the operators $\Cu$ on Hardy spaces and on weighted Dirichlet spaces induce by  complex Borel measures on $[0, 1)$.  Galanopoulos,  Girela et al. \cite{03} studied  the behaviour of the operators $\Cu$  on the Dirichlet space and on the
analytic Besov spaces. Recently, Sun, Ye et al. \cite{syl} studied the operator $\Cu$ from Besov spaces to $X$, where $X$ is a Banach space of analytic
functions in $\dd$ with $\Lambda^{s}_{\frac{1}{s}}\subseteq X \subseteq \mathcal {B}$.  Bao, Guo et al. \cite{affa}  completely characterized the measures $\mu$ such that $\Cu$ is bounded (compact) on Dirichlet space. In \cite{gxt}, the authors of this paper also considered the boundedness and compactness of $\Cu$ between Bergman space and Bloch space.  Beltr\'{a}n-Meneu, Bonet and Jord\'{a} \cite{bel} systematically investigated the operator $\Cu$ on weighted Banach spaces of analytic function.
The operators  $\Cu$ associated to arbitrary complex Borel measures on  $\dd$  the reader is referred to  \cite{04,05}.

The Bolch type spaces $\mathcal {B}^{\alpha}$ are closely connected to many   analytic function spaces, such as Bergman space, Korenblum space, Lipschitz space, $F(p,q,s)$ space, mixed norm space et al.  Therefore, the operator $\Cu$  acting between  Bloch type spaces can serve as a good model when we study  the  operator $\Cu$ on the spaces.
In this paper we  study the action of  the operator $\Cu$ between Bloch type spaces. The operator $\Cu$  on  such  spaces do not seem to have been studied extensively in the past, so we have attempted to collect here the consequences of applying to them various standard techniques of analysis.   

 The Carleson-type measures play a basic role in the studies of  $\Cu$. Let $I\subset \partial \mathbb{D}$ be an arc, and  $\vert I\vert $  denote the length of $I$. The Carleson square $S(I)$ is defined as
$$S(I)=\{re^{i\vartheta}:e^{i\vartheta}\in I,\ 1-\frac{\vert I\vert }{2\pi}\leq r<1\}.$$

Let $\mu$ be a positive Borel measure on $\mathbb{D}$. For $0\leq \beta<\infty$ and $0<t<\infty$, we
say that $\mu$ is a $\beta$-logarithmic $t$-Carleson measure (resp. a vanishing $\beta$-logarthmic $t$-Carleson measure) if
$$
\sup_{ I \subset \partial\mathbb{D}}\frac{\mu(S(I))(\log\frac{2\pi}{\vert I\vert })^{\beta}}{\vert I\vert ^{t}}<\infty,\ \ \left(\mbox{resp.}\ \ \lim_{\vert I\vert \rightarrow0}\frac{\mu(S(I))(\log\frac{2\pi}{\vert I\vert })^{\beta}}{\vert I\vert ^{t}}=0\right).
$$
If $\beta=0$ and $t=1$, we  say  that $\mu$ is a Carleson measure. See \cite{H16} for more about  logarithmic type Carleson measure.

A positive Borel measure $\mu$ on $[0,1)$ can be seen as a Borel measure on $\mathbb{D}$ by identifying it with the measure $\overline{\mu}$ defined by
$$
 \overline{\mu}(E)=\mu(E\cap [0,1)), \ \ \mbox{for any Borel subset }\ E \ \ \mbox{of}\  \ \mathbb{D}.
$$

In this way, a positive Borel measure $\mu$ on $[0,1)$ is a $\beta$-logarithmic $t$-Carleson measure if and only if there exists a constant $M>0$ such that
$$
\mu([s,1))\log^{\beta}\frac{e}{1-s}\leq M(1-s)^{t},\ \ 0\leq s<1.
$$


 Throughout the paper, the letter $C$ will denote an absolute constant whose value depends on the parameters
indicated in the parenthesis, and may change from one occurrence to another. We will use
the notation $``P\lesssim Q"$ if there exists a constant $C=C(\cdot) $ such that $`` P \leq CQ"$, and $`` P \gtrsim Q"$ is
understood in an analogous manner. In particular, if  $``P\lesssim Q"$  and $ ``P \gtrsim Q"$ , then we will write $``P\asymp Q"$.

 \section{Preliminaries} \label{prelim}

In this section, we present  some preliminary results needed for the rest of the paper. We start with the following lemma which can be found, for example, in \cite{b7}.
 \begin{lemma}\label{lem2.1}
Let  $0<\alpha<\infty$  and $f\in \ba$. Then for each $z\in \dd$, we have the following
inequalities:
\begin{equation*}
|f(z)|\lesssim
\begin{cases}
||f||_{\ba},  \text{ if  $0<\alpha<1;$}\\
||f||_{\ba}\log\frac{2}{1-|z|},  \text{ if  $\alpha=1;$ }\\
\frac{||f||_{\ba}}{(1-|z|)^{\alpha-1}}, \text{ if $\alpha>1.$}
\end{cases}
\end{equation*}
\end{lemma}

The following result follows from Corollary 3.2 in \cite{IB} or Theorem 2.26 in \cite{Wulan1}.
\begin{lemma}\label{lm2.2}
Let  $\alpha>0$ and $f\in H(\mathbb{D})$, $f(z)=\sn a_{n}z^{n}$, $a_{n}\geq 0 $ for all $n\geq 0$. Then $f\in \mathcal {B}^{\alpha}$ if and only if
$$\sup_{n\geq 1}n^{-\alpha}\sum_{k=1}^{n}ka_{k}<\infty.$$
\end{lemma}


 The following result follows from Theorem 2.1  and Theorem 2.4 in \cite{H10}.
\begin{lemma}
Let $0<s<\infty$ and $\mu$ be a finite positive Borel measure on the interval $[0,1)$. Then the following statements hold:
\\ (1)\ $\mu$ is an $s$-Carleson measure if and only if $\mu_{n}= O(\frac{1}{n^{s}})$.
\\ (2)\  $\mu$ is a vanishing $s$-Carleson measure if and only if $\mu_{n}= o(\frac{1}{n^{s}})$.
\end{lemma}

The following integral estimates are  practical. Although we only use a partial case in this article, we present a complete result here for the reader's reference.
\begin{lemma}
Let $\delta>-1$, $c\geq 0$ and  $k$ be a real number. Then the integral
$$I_{r}=\ii \frac{(1-t)^{\delta}}{(1-tr)^{\delta+c+1}}\log^{k}\frac{e}{1-t}dt, \ \ (0 \leq r<1)$$
have the following asymptotic properties:
\\ (1)\ If $c=0$ and $k<-1$, then $I_{r}\asymp 1 $;
\\ (2)\ If $c=0$ and $k=-1$, then $I_{r}\asymp \log\log\frac{e^{2}}{1-r} $;
\\ (3)\ If $c=0$ and $ k>-1$, then $I_{r}\asymp \log^{k+1}\frac{e}{1-r}$;
\\ (4)\ If $c>0$, then  $I_{r}\asymp \frac{1}{(1-r)^{c}}\log^{k}\frac{e}{1-r}$.
\end{lemma}
\begin{proof}
The proof of  (3)-(4) is stated in  \cite[Lemma 2.2 ]{ces17}. We just need to consider the case $c=0$ and $k\leq -1$.

Without loss of generality, we may assume that $1-\frac{e}{8}<r<1$. Let $x=\frac{r(1-t)}{1-tr}$, then
 \[ \begin{split}
\ii \frac{(1-t)^{\delta}}{(1-tr)^{\delta+1}}\log^{k}\frac{e}{1-t}dt&=\int_{0}^{r}\frac{x^{\delta}}{r^{\delta+1}(1-x)}\log^{k}\frac{er(1-x)}{(1-r)x}dx\\
& \asymp \int_{0}^{\frac{1}{2}}x^{\delta}\log^{k}\frac{e}{(1-r)x}dx+ \int_{\frac{1}{2}}^{r}\frac{1}{1-x}\log^{k}\frac{e(1-x)}{1-r}dx\\
&= \frac{1}{(1-r)^{\delta+1}}\int_{0}^{\frac{1-r}{2}}y^{\delta}\log^{k}\frac{e}{y}dy+ \int_{2(1-r)}^{1}\frac{1}{y}\log^{k}\frac{e}{y}dy.
  \end{split} \]
It is clear that
$$\lim_{r\rightarrow 1^{-}}\frac{\int_{0}^{\frac{1-r}{2}}y^{\delta}\log^{k}\frac{e}{y}dy}{(1-r)^{\delta+1}\log^{k}\frac{e}{1-r}}= \frac{1}{(\delta+1)2^{\delta+1}}.$$
This implies that
$$ \frac{1}{(1-r)^{\delta+1}}\int_{0}^{\frac{1-r}{2}}y^{\delta}\log^{k}\frac{e}{y}dy \asymp  \log^{k}\frac{e}{1-r} \ (r\rightarrow 1^{-}). \eqno{(2.1)}$$
At the same time,
$$ \int_{2(1-r)}^{1}\frac{1}{y}\log^{-1}\frac{e}{y}dy=\log\log\frac{e}{2(1-r)}\asymp \log\log\frac{e^{2}}{1-r}\   (r\rightarrow 1^{-}). \eqno{(2.2)}$$
When $k<-1$, we have
$$ \int_{2(1-r)}^{1}\frac{1}{y}\log^{k}\frac{e}{y}dy\leq \ii \frac{1}{y}\log^{k}\frac{e}{y}dy=\frac{-1}{k+1}. \eqno{(2.3)}$$
By (2.1)-(2.3) we may obtain that (1) and (2) hold.
\end{proof}

The following lemma is a direct consequence of Theorem 3.1 in \cite{lin}.
\begin{lemma}
 Let $0<\alpha, \beta<\infty$. Suppose $T$ is a bounded operator from $\ba$ into $\bb$, then $T$  is a compact operator from $\ba$ into $\bb$ if and only if  for any bounded sequence $\{h_{n}\}$ in $\ba$ which converges to $0$ uniformly on every compact subset of $\mathbb{D}$, we have
$\displaystyle{\lim_{n\rightarrow \infty}||T(h_{n})|| _{\bb}=0}$.
\end{lemma}

\section{The boundedness of  $\Cu$  acting between Bloch type spaces} \label{sec:HL-spaces}

We now study the  boundedness of  $\Cu$  acting between Bloch type spaces.

\begin{theorem}
Let $\mu$ be a  finite positive Borel measure on the interval $[0,1)$. If $0<\alpha<1$ and $0<\beta<2$, then    the following  conditions are  equivalent.
\\(1)  $\Cu: \mathcal {B}^{\alpha}\rightarrow \mathcal {B}^{\beta}$ is bounded.
\\(2) $\Cu: \mathcal {B}^{\alpha}\rightarrow \mathcal {B}^{\beta}$ is compact.
\\(3)  The measure $\mu$ is a $2-\beta$ Carleson measure.
\end{theorem}
\begin{proof}
The implication  of $ (2)\Rightarrow(1)$ is obvious.

$ (1)\Rightarrow(3)$.  Suppose  $\Cu: \mathcal {B}^{\alpha}\rightarrow \mathcal {B}^{\beta}$ is bounded. Let $f(z)=\sum_{n=1}^{\infty}n^{\alpha-2}z^{n}$, then
$$|f'(z)|= \left|\sum_{n=0}^{\infty}(n+1)^{\alpha-1}z^{n}\right|\leq \sum_{n=0}^{\infty}(n+1)^{\alpha-1}|z|^{n} \lesssim \frac{1}{(1-|z|)^{\alpha}}.$$
This means that  $f\in \ba$. Since
 $$\Cu(f)(z)=\sum_{n=1}^{\infty} \mu_{n}\left(\sum_{k=1}^{n}k^{\alpha-2}\right)z^{n}\in \bb$$
 and the sequence  $\left\{ \mu_{n}\left(\sum_{k=1}^{n}k^{\alpha-2}\right)\right\}^{\infty}_{n=1}$ is a  nonnegative sequence,  it follows from Lemma \ref{lm2.2} that
$$\sup_{n\geq 1 } n^{-\beta}\sum_{k=1}^{n}k\mu_{k}\left(\sum_{j=1}^{k}j^{\alpha-2}\right)<\infty.$$
Since $\alpha \in (0,1)$, for each $n\geq 1$, it follows that
\[ \begin{split}
1 &\gtrsim n^{-\beta}\sum_{k=1}^{n}k\mu_{k}\left(\sum_{j=1}^{k}j^{\alpha-2}\right)\\
& \gtrsim \mu_{n}n^{-\beta}\sum_{k=1}^{n}k \asymp \mu_{n}n^{2-\beta}.
\end{split} \]
 Lemma 2.3 shows that  $\mu$ is a $2-\beta$ Carleson measure.

$(3)\Rightarrow(1)$.  Assume $\mu$ is a $2-\beta$ Carleson measure.  Since $0<\alpha<1$, the integral $\ii \frac{dt}{(1-t)^{\alpha}}<\infty$.  Thus, for any given $\varepsilon>0$, there exists $0<t_{0}<1$ such that
 $$ \int_{t_{0}}^{1}\frac{dt}{(1-t)^{\alpha}}<\varepsilon .\eqno{(3.1)}$$
This also yields  that
 $$1-t_{0}<\varepsilon^{\frac{1}{1-\alpha}}. \eqno{(3.2)}$$
Let $\{f_{n}\}_{n=1}^{\infty}$ be a bounded
sequence in $\ba$  which converges to $0$ uniformly on every compact subset of  $\mathbb{D}$. Without loss of generality, we may assume  that $\sup_{n\geq 1}||f_{n}||_{\mathcal {B}^{\alpha}}\leq 1$. By the integral representation of  $\Cu$  we see that
\[ \begin{split}
\Cu(f_{n})(z)&=\ii\frac{f_{n}(tz)-f_{n}(t)}{1-tz}d\mu(t)+\ii\frac{f_{n}(t)}{1-tz}d\mu(t)\\
&:=\mathcal {J}_{\mu}(f_{n})(z)+ \mathcal {I}_{\mu}(f_{n})(z).
\end{split} \]
It follows that
$$||\Cu(f_{n})||_{\mathcal {B}^{\beta}}\leq ||\mathcal {J}_{\mu}(f_{n})||_{\mathcal {B}^{\beta}}+ ||\mathcal {I}_{\mu}(f_{n})||_{\mathcal {B}^{\beta}}.$$
Note  that  the second part of the right-hand side is the integral type Hilbert operator (see e.g. \cite{H14,IB} for the definition).  Therefore,  Corollary 5.4  in \cite{IB} shows that the integral type Hilbert operator $\mathcal {I}_{\mu}$ is compact from  $\mathcal {B}^{\alpha}$ to $ \mathcal {B}^{\beta}$ whenever $0<\alpha<1$ and $0<\beta<2$. This implies that
$$\lim_{n\rightarrow \infty}||\mathcal {I}_{\mu}(f_{n})||_{\mathcal {B}^{\beta}}=0.$$
To complete the proof, it is suffices to prove that
$\lim_{n\rightarrow \infty}||\mathcal {J}_{\mu}(f_{n})||_{\mathcal {B}^{\beta}}=0$ by Lemma 2.5. It is easy to see that
$$
|\mathcal {J}_{\mu}(f_{n})'(z)|\leq   \ii G^{z}_{n}(t)d\mu(t)
$$
where  $$G^{z}_{n}(t)=\frac{|f_{n}'(tz)|}{|1-tz|}+\frac{|f_{n}(tz)-f_{n}(t)|}{|1-tz|^{2}}, \ \  z\in \dd.$$
The Cauchy's integral theorem implies that $\{f'_{n}\}_{n=1}^{\infty}$  converges to $0$ uniformly on every compact subset of  $\mathbb{D}$.
Hence
$$\sup_{z\in\dd}(1-|z|^{2})^{\beta}\int_{0}^{t_{0}}G^{z}_{n}(t)d\mu(t)\lesssim \sup_{|w|\leq t_{0}}(|f'_{n}(w)|+|f_{n}(w)|)\rightarrow 0, \ \ \mbox{as}\ n\rightarrow \infty.$$
Since $\mathcal {B}^{\alpha}\cong\Lambda_{1-\alpha}$,  we have
$$|f_{n}(tz)-f_{n}(t)|\leq t|1-z|^{1-\alpha}||f_{n}||_{\Lambda_{1-\alpha}}\lesssim  |1-z|^{1-\alpha}.$$
For $0<t<1$ and $z\in \dd$, the inequalities
$$\frac{|1-z|}{|1-tz|}\leq \frac{1-t}{|1-tz|}+\frac{|t-z|}{|1-tz|}\leq 2 $$ imply that
$$\frac{|f_{n}(tz)-f_{n}(t)|}{|1-tz|^{2}} \lesssim \frac{ 1}{|1-tz|^{1+\alpha}}.\eqno{(3.3)}$$
By the definition of $\ba$ and (3.3) we have
\[ \begin{split}
\int_{t_{0}}^{1}G^{z}_{n}(t)d\mu(t) & \lesssim  \int_{t_{0}}^{1}\left (\frac{1}{(1-t|z|)^{\alpha}|1-tz|}+\frac{1}{|1-tz|^{1+\alpha}}\right)d\mu(t)\\
& \lesssim \int_{t_{0}}^{1} \frac{d\mu(t)}{(1-t|z|)^{\alpha+1}}.
\end{split} \]
 Bearing in the mind that  $\mu$ is a $2-\beta$ Carleson measure and that there exists $0<\delta<1$ such that $(1-|z|^{2})^{\beta}<\varepsilon$ for all $\delta<|z|<1$, by   integrating by parts (see \cite[Theorem 5]{ces1}) and using (3.1)-(3.2) we have
\[ \begin{split}
&\ \ \ \ \sup_{z\in\dd}(1-|z|^{2})^{\beta}\int^{1}_{t_{0}}G^{z}_{n}(t)d\mu(t)\\
& \lesssim
\sup_{z\in\dd}(1-|z|^{2})^{\beta}\int_{t_{0}}^{1} \frac{d\mu(t)}{(1-t|z|)^{\alpha+1}}\\
&= \sup_{z\in\dd}(1-|z|^{2})^{\beta}\left(\frac{\mu([t_{0},1))}{(1-t_{0}|z|)^{\alpha+1}}+(\alpha+1)|z|\int_{t_{0}}^{1}\frac{\mu([t,1))}{(1-t|z|)^{\alpha+2}}dt\right)\\
&  \lesssim  \left(\sup_{|z|\leq \delta}+\sup_{\delta<|z|<1}\right)(1-|z|^{2})^{\beta}\frac{\mu([t_{0},1))}{(1-t_{0}|z|)^{\alpha+1}} +\sup_{z\in\dd}\int_{t_{0}}^{1}\frac{(1-t)^{2-\beta}(1-|z|^{2})^{\beta}}{(1-t|z|)^{\alpha+2}}dt\\
&  \lesssim  (1-t_{0})^{2-\beta}+\varepsilon+\int_{t_{0}}^{1}\frac{dt}{(1-t)^{\alpha}}\\
& \lesssim \varepsilon^{\frac{2-\beta}{1-\alpha}} +\varepsilon .
\end{split} \]
Consequently,
$$\lim_{n\rightarrow \infty}||\mathcal {J}_{\mu}(f_{n})||_{\mathcal {B}^{\beta}}=0.$$
This implies that
$\Cu: \mathcal {B}^{\alpha}\rightarrow \mathcal {B}^{\beta}$ is compact.
\end{proof}

\begin{theorem}
Let $\mu$ be a  finite positive Borel measure on the interval $[0,1)$. If $\alpha>1$ and $0<\beta<\alpha+1$, then $\Cu: \mathcal {B}^{\alpha}\rightarrow \mathcal {B}^{\beta}$ is bounded if and only if $\mu$ is an $\alpha+1-\beta$ Carleson measure.
\end{theorem}
\begin{proof}
Suppose  $\Cu: \mathcal {B}^{\alpha}\rightarrow \mathcal {B}^{\beta}$ is bounded.
For $0<a<1$, let
$$f_{a}(z)=\frac{(1-a)}{(1-az)^{\alpha}}.$$
Then it is easy to check that $\sup_{0<a<1}||f_{a}||_{\ba}\lesssim 1$. By the integral form of $\Cu$ we get
$$\Cu(f_{a})'(z)=\ii\frac{tf'_{a}(tz)}{1-tz}d\mu(t)+\ii\frac{tf_{a}(tz)}{(1-tz)^{2}}d\mu(t).$$
The boundedness of $\Cu$  and Lemma 2.1 imply that
 $$|\Cu(f_{a})'(z)|\lesssim \frac{||f_{a}||_{\ba}}{(1-|z|)^{\beta}}.$$
Therefore, for any $\frac{1}{2}<a<1$ we have
\[ \begin{split}
 \frac{1}{(1-a)^{\beta}}&\gtrsim|\Cu(f_{a})'(a)|
= \Cu(f_{a})'(a)\\
&\geq \ii \frac{tf_{a}(ta)}{(1-ta)^{2}}d\mu(t)\\
&= (1-a)\ii \frac{td\mu(t)}{(1-ta)^{2}(1-a^{2}t)^{\alpha}}\\
& \gtrsim  (1-a)\int_{a}^{1} \frac{d\mu(t)}{(1-ta)^{2}(1-a^{2}t)^{\alpha}}\\
& \gtrsim \frac{\mu([a,1))}{(1-a)^{\alpha+1}}.
  \end{split} \]
This gives that
$$\mu([a,1))\lesssim (1-a)^{\alpha+1-\beta} \ \ \mbox{for all}\ \ \frac{1}{2}<a<1. $$
Hence $\mu$ is an $\alpha+1-\beta$ Carleson measure.
\comment{
Let $f(z)=\frac{1}{(1-z)^{\alpha-1}}$, it is easy to check that $f\in \mathcal {B}^{\alpha}$ and
$$\Cu(f)(z)=\sum_{n=0}^{\infty}a_{n}z^{n},\ \ z\in \dd,$$
where
$$a_{n}=\mu_{n}\left(\sum_{k=0}^{n}\frac{\Gamma(k+\alpha-1)}{\Gamma(k+1)\Gamma(\alpha-1)}\right).$$
It is obvious that $\{a_{n}\}$ is  a nonnegative sequence. Since $\Cu(f)\in \bb$, it follows from Lemma \ref{lm2.2} that
$$\sup_{n\geq 1}n^{-\beta}\sum_{k=1}^{n}ka_{k}<\infty.$$
By the Stirling formula,
$$
\frac{\Gamma(n+\alpha-1)}{\Gamma(\alpha-1)\Gamma(n+1)}\asymp (n+1)^{\alpha-2}
$$
for all nonnegative integers $n$. This together with $\{\mu_{n}\}$ is decreasing with $n$  we deduce that
\[ \begin{split}
n^{-\beta}\sum_{k=1}^{n}ka_{k}
&=n^{-\beta}\sum_{k=1}^{n}k\mu_{k}\left(\sum_{j=0}^{k}\frac{\Gamma(j+\alpha-1)}{\Gamma(j+1)\Gamma(\alpha-1)}\right)\\
& \asymp n^{-\beta}\sum_{k=1}^{n}k\mu_{k}\left(\sum_{j=0}^{k} (j+1)^{\alpha-2}\right)\\
&\gtrsim  n^{-\beta}\mu_{n}\sum_{k=1}^{n} k^{\alpha}\asymp \mu_{n}n^{\alpha+1-\beta}.
  \end{split} \]
This implies that  $\mu_{n}= O\left(\frac{1}{n^{\alpha+1-\beta}}\right)$ for all $n\geq 1$. Lemma 2.3 shows that  $\mu$ is an $\alpha+1-\beta$ Carleson measure.
}

Conversely, suppose $\mu$ is an $\alpha+1-\beta$ Carleson measure and $f\in \ba$. Using the   integral form of $\Cu(f)$ and  Lemma  \ref{lem2.1} we deduce that
$$
|\Cu(f)'(z)|\leq  \ii\frac{|f'(tz)|}{|1-tz|}d\mu(t)+\ii\frac{|f(tz)|}{|1-tz|^{2}}d\mu(t)$$
$$ \lesssim||f||_{\ba} \ii\frac{d\mu(t)}{(1-t|z|)^{\alpha+1}}. \ \ \ \ \eqno{(3.4)}
$$
Take $z\in \dd$ and let $|z|=r$.
Integrating by parts and using the fact that $\mu$ is an  $\alpha+1-\beta$   Carleson measure and Lemma 2.4, we obtain
\[ \begin{split}
 \ii\frac{d\mu(t)}{(1-tr)^{\alpha+1}}&= \mu([0,1))+(\alpha+1)r\ii\frac{\mu([t,1))}{(1-tr)^{\alpha+2}}dt\\
 & \lesssim \mu([0,1))+\ii\frac{(1-t)^{\alpha+1-\beta}}{(1-tr)^{\alpha+2}}dt\\
 & \lesssim\mu([0,1)) + \frac{1}{(1-r)^{\beta}}.
  \end{split} \]
This together with  (3.4) imply that  $\Cu: \mathcal {B}^{\alpha}\rightarrow \mathcal {B}^{\beta}$ is bounded.
\end{proof}
\begin{theorem}
Let $\mu$ be a  finite positive Borel measure on the interval $[0,1)$. If $0<\beta\leq2$, then  $\Cu: \mathcal {B}\rightarrow \mathcal {B}^{\beta}$ is bounded if and only if
$$\sup_{0<t<1}\frac{\mu([t,1))\log\frac{e}{1-t}}{(1-t)^{2-\beta}}<\infty.\eqno{(3.5)}$$
\end{theorem}
\begin{proof}
Suppose  $\Cu: \mathcal {B}\rightarrow \mathcal {B}^{\beta}$ is bounded. Let $f(z)=\log\frac{1}{1-z}$, it is clear that $f\in \mathcal {B}$ and
$$\Cu(f)(z)=\sum_{n=1}^{\infty}\mu_{n}\left(\sum_{k=1}^{n}\frac{1}{k}\right)z^{n},\ \ z\in \dd,$$
Since $\Cu(f)\in \bb$, by the definition of $\bb$ we have that
$$\sum_{n=1}^{\infty}n\mu_{n}\left(\sum_{k=1}^{n}\frac{1}{k}\right)r^{n-1}\lesssim \frac{1}{(1-r)^{\beta}}, \ \ 0<r<1.$$
For $N\geq 2$ take $r_{N}=1-\frac{1}{N}$. Since the sequence $\{\mu_{k}\}$ is decreasing,
simple estimations lead us to the following
\[ \begin{split}
N^{\beta}&\gtrsim \sum_{n=1}^{\infty}n\mu_{n}\left(\sum_{k=1}^{n}\frac{1}{k}\right)r_{N}^{n-1}\\
&\gtrsim  \sum_{n=1}^{N}n\mu_{n}\left(\sum_{k=1}^{n}\frac{1}{k}\right)r_{N}^{N-1}\\
&\gtrsim \mu_{N} \sum_{n=1}^{N}n\log(n+1)\\
&\gtrsim \mu_{N}N^{2}\log(N+1)
  \end{split} \]
  This implies that $\mu_{N}=O\left(\frac{1}{N^{2-\beta}\log(N+1)}\right)$.
  \comment{
$$\sup_{n\geq 1}n^{-\beta}\sum_{k=1}^{n}kb_{k}<\infty.$$
Simple estimates give that
\[ \begin{split}
n^{-\beta}\sum_{k=1}^{n}kb_{k}
&=n^{-\beta}\sum_{k=1}^{n}k\mu_{k}\left(\sum_{j=1}^{k}\frac{1}{j}\right)\\
& \asymp n^{-\beta}\sum_{k=1}^{n}k\mu_{k}\log(k+1)\\
&\gtrsim  n^{-\beta}\mu_{n}\sum_{k=1}^{n} k\log(k+1)\asymp \mu_{n}n^{2-\beta}\log(n+1).
  \end{split} \]
Hence $\mu_{n}=O(\frac{1}{n^{2-\beta}\log(n+1)})$ for all $n\geq 1$.
}
 The desired result follows from the  inequalities
$$\mu([1-\frac{1}{N},1))\lesssim\int^{1}_{1-\frac{1}{N}}t^{N}d\mu(t)\leq \ii t^{N}d\mu(t)\lesssim \frac{1}{N^{2-\beta}\log(N+1)}.$$

Conversely, suppose (3.5) holds.
Integrating by parts  we have
\[ \begin{split}
\ii t^{n}d\mu(t)&=n\ii t^{n-1} \mu([t,1))dt\\
& \lesssim n\ii t^{n-1}(1-t)^{2-\beta}\log^{-1}\frac{e}{1-t}dt.\\
  \end{split} \]
  Let $\phi(t)=(1-t)^{2-\beta}\log^{-1}\frac{e}{1-t}$, then  $\phi(t)$ is  regular in the sense of Pel\'{a}ez and R\"{a}tty\"{a} \cite{pl}. Then, using Lemma 1.3 and (1.1) in \cite{pl}, we have
 $$n\ii t^{n-1}\phi(t)dt\asymp \phi\left(1-\frac{1}{n}\right)\asymp \frac{1}{n^{2-\beta}\log(n+1)}.$$
This implies that
$$\mu_{n}\lesssim \frac{1}{n^{2-\beta}\log(n+1)}\ \  \mbox{for all}\ \ n \geq 1 .\eqno{(3.6)}$$
Let  $f(z)=\sn a_{n}z^{n}\in \mathcal {B}$, it follows from  Corollary D in \cite{Ka} that
$$\left|\sum^{n}_{k=1}a_{k}\right|\lesssim ||f||_{\mathcal {B}}\log(n+1).$$
By (3.6) and above inequality we get
\[ \begin{split}
 (1-|z|^{2})^{\beta}|\Cu(f)'(z)|
&= (1-|z|^{2})^{\beta}\left|\sum_{n=1}^{\infty}n\mu_{n}\left(\sum^{n}_{k=1}a_{k}\right)z^{n-1}\right|\\
& \lesssim (1-|z|^{2})^{\beta} \sum_{n=1}^{\infty}n\mu_{n}\left|\sum^{n}_{k=1}a_{k}\right||z|^{n-1}\\
& \lesssim ||f||_{\mathcal {B}} (1-|z|^{2})^{\beta} \sum_{n=1}^{\infty}n\mu_{n}\log(n+1)|z|^{n-1}\\
&\lesssim ||f||_{\mathcal {B}} (1-|z|^{2})^{\beta}  \sum_{n=1}^{\infty}n^{\beta-1}|z|^{n-1}\\
&\lesssim||f||_{\mathcal {B}}.
  \end{split} \]
 This shows that
  $$\sup_{z\in \dd}(1-|z|^{2})^{\beta} |\Cu(f)'(z)|\lesssim ||f||_{\mathcal {B}}.$$
Hence, we have that  $\Cu: \mathcal {B}\rightarrow \mathcal {B}^{\beta}$ is bounded.
\comment{
By Lemma 2.1 and (3.1) we have
\[ \begin{split}
|\Cu(f)'(z)| & \leq  \ii\frac{|f'(tz)|}{|1-tz|}d\mu(t)+\ii\frac{|f(tz)|}{|1-tz|^{2}}d\mu(t)\\
& \lesssim ||f||_{\mathcal {B}}\left(\ii\frac{d\mu(t)}{(1-t|z|)^{2}}+\ii\frac{\log\frac{e}{1-t|z|}}{(1-t|z|)^{2}}d\mu(t)\right)\\
& \lesssim ||f||_{\mathcal {B}} \log\frac{e}{1-|z|}\ii\frac{d\mu(t)}{(1-t|z|)^{2}}.
  \end{split} \]
Integrating by parts and using the fact that  $\mu$  a  $1$-logarithmic $2-\beta$ Carleson measure and
Lemma 2.4, we obtain
\[ \begin{split}
\ii\frac{d\mu(t)}{(1-t|z|)^{2}} &=  \mu([0,1))+2|z|\ii\frac{\mu([t,1))}{(1-t|z|)^{3}}dt\\
&  \lesssim \mu([0,1)) + \ii\frac{(1-t)^{2-\beta}}{(1-t|z|)^{3}}\log^{-1}\frac{e}{1-t}dt\\
& \lesssim \mu([0,1)) + \frac{ \log^{-1}\frac{e}{1-|z|}}{(1-|z|)^{\beta}}.
  \end{split} \]
  }
\end{proof}

\comment{
\begin{theorem}
Let $\mu$ be a  finite positive Borel measure on the interval $[0,1)$. Then $\Cu: \mathcal {B}\rightarrow \mathcal {B}^{2}$ is bounded if and only if
$$\mu([t,1))\lesssim \log^{-1}\frac{2}{1-t}\ \mbox{for all}\ 0\leq t<1.$$
\end{theorem}
\begin{proof}
 Suppose $\Cu: \mathcal {B}\rightarrow \mathcal {B}^{2}$ is bounded, then arguing as Theorem 3.3  may obtain
 $$\mu_{n}=O\left(\frac{1}{\log(n+1)}\right)\ \ \mbox{for all}\ n\geq 1.$$
 This implies that $$\mu([t,1))\lesssim \log^{-1}\frac{2}{1-t}\ \mbox{for all}\ 0\leq t<1.$$
The proof of sufficiency is similar to Theorem 3.3 and hence we omit  the details.

\end{proof}
}

There are  three cases left: (i) $\alpha=1$ and $\beta>2$; (ii) $0<\alpha<1$ and $\beta\geq 2$;  (iii) $\alpha>1$ and $\beta\geq \alpha +1$. We show that the operator $\Cu$   is always  a bounded operator from $\ba$ to $\bb$  in these cases.
\begin{theorem}
Let $\mu$ be a  finite positive Borel measure on the interval $[0,1)$. If $\alpha$ and $\beta$ satisfies one of the conditions (i)-(iii), then $\Cu: \mathcal {B}^{\alpha}\rightarrow \mathcal {B}^{\beta}$  is  bounded.
\end{theorem}
\begin{proof}
We  only prove the case of (i), since the proofs  of other cases are similar. Let $f\in \mathcal {B}$, then
\[ \begin{split}
|\Cu(f)'(z)| & \leq  \ii\frac{|f'(tz)|}{|1-tz|}d\mu(t)+\ii\frac{|f(tz)|}{|1-tz|^{2}}d\mu(t)\\
& \lesssim ||f||_{\mathcal {B}} \ii \frac{d\mu(t)}{(1-t|z|)|1-tz|}+ ||f||_{\mathcal {B}} \ii \frac{\log\frac{e}{1-t|z|}}{|1-tz|^{2}}d\mu(t)\\
& \lesssim ||f||_{\mathcal {B}}  \ii \frac{\log\frac{e}{1-t|z|}}{(1-t|z|)^{2}}d\mu(t).
  \end{split} \]
If $\beta>2$, it is clear that
$$\sup_{z\in \dd}(1-|z|)^{\beta-2}\log\frac{e}{1-|z|}<\infty.$$
This implies that
$$\sup_{z\in \dd}(1-|z|^{2})^{\beta}|\Cu(f)'(z)| \lesssim ||f||_{\mathcal {B}}.$$
\comment{ We  only prove the case of (i), since the proof  of other cases are similar. Let $f\in \ba$, if $0<\alpha<1$, then $f\in H^{\infty}$. Therefore,
\[ \begin{split}
|\Cu(f)'(z)| & \leq  \ii\frac{|f'(tz)|}{|1-tz|}d\mu(t)+\ii\frac{|f(tz)|}{|1-tz|^{2}}d\mu(t)\\
& \lesssim ||f||_{\ba} \ii \frac{d\mu(t)}{(1-t|z|)^{\alpha+1}}+ ||f||_{\ba} \ii \frac{d\mu(t)}{(1-t|z|)^{2}}\\
& \lesssim ||f||_{\ba}  \ii \frac{d\mu(t)}{(1-t|z|)^{2}}\lesssim \frac{  ||f||_{\ba} }{(1-|z|)^{2}}.
  \end{split} \]
If $\beta\geq 2$, then $$\sup_{z\in \dd}(1-|z|^{2})^{\beta}|\Cu(f)'(z)| \lesssim  ||f||_{\ba}.  $$
This means that $\Cu: \mathcal {B}^{\alpha}\rightarrow \mathcal {B}^{\beta}$  is  bounded.}
The proof is complete.\end{proof}

Let  $1\leq p<\infty$ and $0<\alpha\leq 1$, the mean Lipschitz space $\Lambda^p_\alpha$ consists of the functions $f\in \hd$ having a non-tangential limit  almost everywhere for which  $\omega_p(t, f)=O(t^\alpha)$ as $t\rightarrow 0$. Here $\omega_p(\cdot, f)$ is the integral modulus of continuity of order $p$ of the function $f(e^{i\theta})$. It is known that  (see e.g., \cite[Chapter 5]{b1}) $\Lambda^p_\alpha$ is a subset of $H^p$ and $\Lambda^p_\alpha$ consists of those functions $f\in \hd$ satisfying
$$
\|f\|_{\Lambda^p_\alpha}=|f(0)|+\sup_{0<r<1}(1-r)^{1-\alpha}M_p(r, f')<\infty.
$$

 Theorem 4.2 in \cite{baoo} and Theorem 3.1  lead to the following corollary.
\begin{corollary}
Let $\mu$ be a  finite positive Borel measure on the interval $[0,1)$, $0<\alpha<1$ and $1<p<\infty$. Let $X$ and $Y$ be two Banach subspaces of $\hd$ with $\ba \subset X\subset H^{\infty}$ and $\Lambda^{p}_{\frac{1}{p}}\subset Y \subset \mathcal {B}$. Then the following conditions are equivalent.
\\(1)\ The measure $\mu$ is a Carleson measure.
\\(2)\ The operator $\Cu$ is bounded from $X$ into $Y$.
\end{corollary}
\begin{proof}
$(1)\Rightarrow (2)$. Assume that $\mu$ is a Carleson measure and take $f\in X$. Since $X\subset H^{\infty}$, we have that  $f\in H^{\infty}$.  Theorem 4.2 in \cite{baoo}  shows that $\Cu (H^{\infty})\subset Y $  if and only if  $\mu$ is a Carleson measure.
 This implies that $\Cu (f)\in Y$.

$(2)\Rightarrow (1)$. Suppose that $\Cu$ is bounded from $X$ into $Y$. Then $\Cu$ is bounded from $\ba$ into $\mathcal {B}$.  Now, Theorem 3.1 shows that $\mu$ is a Carleson measure.
\end{proof}

\section{ The compactness of  $\Cu$  acting between Bloch type  spaces}

\begin{theorem}
Let $\mu$ be a  finite positive Borel measure on the interval $[0,1)$. If $\alpha>1$ and $0<\beta<\alpha+1$, then $\Cu: \mathcal {B}^{\alpha}\rightarrow \mathcal {B}^{\beta}$ is compact if and only if $\mu$ is a vanishing $\alpha+1-\beta$ Carleson measure.
\end{theorem}
\begin{proof}
Assume that $\Cu: \mathcal {B}^{\alpha}\rightarrow \mathcal {B}^{\beta}$ is compact. For $0<a<1$, set
$$f_{a}(z)=\frac{1-a}{(1-az)^{\alpha}}, \ z\in \dd.$$
Then it is clear that $f_{a}\in \ba$ for all $0<a<1$ and $\sup_{0<a<1}||f_{a}||_{\ba}\lesssim 1$.
Moreover, $f_{a}\rightarrow 0$, as $a\rightarrow 1$, uniformly on compact subset of $\dd$.  Lemma 2.5 implies that
$$||\Cu(f_{a})||_{\bb}\rightarrow 0, \ \ \mbox{as}\ a\rightarrow 1. \eqno{(4.1)}$$
\comment{
For $0<a<1$ and $z\in \dd$, it follows from (3.1) that
$$\Cu(f_{a})'(z)=\ii\frac{tf'_{a}(tz)}{1-tz}d\mu(t)+\ii\frac{tf_{a}(tz)}{(1-tz)^{2}}d\mu(t).$$
Therefore,
\[ \begin{split}
|\Cu(f_{a})'(a)| &= \Cu(f_{a})'(a)\geq \ii \frac{tf_{a}(ta)}{(1-ta)^{2}}d\mu(t)\\
&= (1-a)\ii \frac{d\mu(t)}{(1-ta)^{2}(1-a^{2}t)^{\alpha}}\\
& \geq  (1-a)\int_{a}^{1} \frac{d\mu(t)}{(1-ta)^{2}(1-a^{2}t)^{\alpha}}\\
& \gtrsim \frac{\mu([a,1))}{(1-a)^{\alpha+1}}.
  \end{split} \]
  }
  Arguing as the proof of Theorem 3.2 we have
$$\mu([a,1))\lesssim (1-a)^{\alpha+1-\beta}||\Cu(f_{a})||_{\bb}.$$
This and (4.1) show that $\mu$  is a vanishing $\alpha+1-\beta$ Carleson measure.

On the other hand, suppose $\mu$  is a vanishing $\alpha+1-\beta$ Carleson measure. Then for any $\varepsilon>0$, there exists $0<t_{0}<1$ such that
$$\mu([t,1))<\varepsilon (1-t)^{\alpha+1-\beta}\ \ \mbox{whenever}\ t_{0}\leq t<1.\eqno{(4.2)}$$
Let $\{f_{k}\}_{k=1}^{\infty}$ be a bounded
sequence in $\ba$  which converges to $0$ uniformly on every compact subset of  $\mathbb{D}$.
It is  sufficient to prove that
$$\lim_{k\rightarrow\infty}||\Cu(f_{k})||_{\bb}=0$$
by Lemma 2.5.
By the integral form of $\Cu(f)$  we have that
\[ \begin{split}
\sup_{z\in\dd}(1-|z|^{2})^{\beta}|\Cu(f_{k})'(z)|&\leq \sup_{z\in\dd}(1-|z|^{2})^{\beta}\ii \frac{|f'_{k}(tz)|}{|1-tz|}d\mu(t)\\
& + \sup_{z\in\dd}(1-|z|^{2})^{\beta}\ii \frac{|f_{k}(tz)|}{|1-tz|^{2}}d\mu(t).
  \end{split} \]
The Cauchy's integral theorem implies that $\{f'_{k}\}_{k=1}^{\infty}$  converges to $0$ uniformly on every compact subset of  $\mathbb{D}$.
This gives
$$\sup_{z\in\dd}(1-|z|^{2})^{\beta}\int_{0}^{t_{0}}\frac{|f'_{k}(tz)|}{|1-tz|}d\mu(t)\lesssim \sup_{|w|\leq t_{0}}|f'_{k}(w)|\rightarrow 0, \ \ \mbox{as}\ k\rightarrow \infty.$$
Note  that there exists $0<\delta<1$ such that $(1-|z|^{2})^{\beta}<\varepsilon$ for all $\delta<|z|<1$, by integrating by parts and using (4.2) and Lemma 2.4 we have
\[ \begin{split}
&\ \ \ \ \sup_{z\in\dd}(1-|z|^{2})^{\beta}\int_{t_{0}}^{1}\frac{|f'_{k}(tz)|}{|1-tz|}d\mu(t)\\
&\lesssim \sup_{z\in\dd}(1-|z|^{2})^{\beta}\int_{t_{0}}^{1}\frac{d\mu(t)}{(1-t|z|)^{\alpha+1}}\\
& = \sup_{z\in\dd}(1-|z|^{2})^{\beta}\left(\frac{\mu([t_{0},1))}{(1-t_{0}|z|)^{\alpha+1}}+(\alpha+1)|z|\int_{t_{0}}^{1}\frac{\mu([t,1))}{(1-t|z|)^{\alpha+2}}dt\right)\\
& \lesssim \left(\sup_{|z|\leq \delta} +\sup_{\delta<|z|<1}\right)(1-|z|^{2})^{\beta}\frac{\mu([t_{0},1))}{(1-t_{0}|z|)^{\alpha+1}}+ \varepsilon \sup_{z\in\dd}(1-|z|^{2})^{\beta}\int_{t_{0}}^{1}\frac{(1-t)^{\alpha+1-\beta}}{(1-t|z|)^{\alpha+2}}dt\\
& \lesssim \varepsilon(1-t_{0})^{\alpha+1-\beta}+\varepsilon+\varepsilon\sup_{z\in\dd}(1-|z|^{2})^{\beta}\ii\frac{(1-t)^{\alpha+1-\beta}}{(1-t|z|)^{\alpha+2}}dt\\
& \lesssim\varepsilon .
  \end{split} \]
Since $\varepsilon$ is arbitrary, it follows that
$$\lim_{k\rightarrow \infty}\sup_{z\in\dd}(1-|z|^{2})^{\beta}\int_{0}^{1}\frac{|f'_{k}(tz)|}{|1-tz|}d\mu(t)=0.$$
Similarly, we can obtain that
$$\lim_{k\rightarrow \infty}\sup_{z\in\dd}(1-|z|^{2})^{\beta}\int_{0}^{1}\frac{|f_{k}(tz)|}{|1-tz|^{2}}d\mu(t)=0.$$
It is obvious that $$\lim_{k\rightarrow\infty}|\Cu(f_{k})(0)|=\mu([0,1))\lim_{k\rightarrow\infty}|f_{k}(0)|=0.$$
Therefore,
$$\lim_{k\rightarrow\infty}||\Cu(f_{k})||_{\bb}=\lim_{k\rightarrow\infty}\left(|\Cu(f_{k})(0)|+\sup_{z\in\dd}(1-|z|^{2})^{\beta}|\Cu(f_{k})'(z)|\right)=0.$$
This means that  $\Cu: \mathcal {B}^{\alpha}\rightarrow \mathcal {B}^{\beta}$ is compact.
\end{proof}
\begin{theorem}
Let $\mu$ be a  finite positive Borel measure on the interval $[0,1)$. If $0<\beta\leq2$, then  $\Cu: \mathcal {B}\rightarrow \mathcal {B}^{\beta}$ is compact if and only if
$$\lim_{t\rightarrow 1}\frac{\mu([t,1))\log\frac{e}{1-t}}{(1-t)^{2-\beta}}=0.$$
\end{theorem}
\begin{proof}
The proof of the sufficiency is similar to  that of Theorem 4.1.  Take the test functions
$$f_{a}(z)=\log^{-1}\frac{2}{1-a}\left(\log\frac{2}{1-az}\right)^{2},\ \  a\in (0,1), \ z\in \dd .$$
Then arguing as the proof of Theorem 4.1 we can obtain the necessity.
\end{proof}



\section*{Conflicts of Interest}
The authors declare that there is no conflict of interest.

\section*{Funding}

 The first author was supported  by the Scientific Research Fund of Hunan Provincial Education Department (NO. 24C0222) and the second author was supported by  the Natural Science Foundation of Hunan Province (No. 2022JJ30369).



\section*{Availability of data and materials}
Data sharing not applicable to this article as no datasets were generated or analysed during
the current study: the article describes entirely theoretical research.

\section*{Acknowledgements}
The authors would like to thank the anonymous referee for
careful review of our paper. The referee has made valuable comments and suggestions
for improving the earlier version of the paper.


\pdfbookmark[1]{ References}{4}

 \end{document}